\documentclass[reqno]{article}
\usepackage{amsmath, amsthm, graphicx}
\theoremstyle{plain}

\usepackage{setspace}
\usepackage{ifthen}
\def\defas{\; \buildrel\rm def\over= \;}

\begin{document}
\bibliographystyle{amsalpha}
\title{A Riemann-Farey Computation}
\author{Scott B. Guthery\\sbg@acw.com}
\maketitle

The Riemann hypothesis is true if and only if
\begin{equation}\label{EQ1}
R(m) = \sum_{i=2}^n\left({F_m(i) - \frac{i}{n}}\right)^2 =
O(m^{-1+\epsilon})
\end{equation}
where $F_m(i)$ is the $i^{th}$ element in the Farey sequence of
order $m$ and
$$
n=\sum_{k=2}^m \phi(k)\text{.}
$$

\pagebreak

\begin{figure}[h!]
\begin{center}
\leavevmode
  \includegraphics[width=100mm]{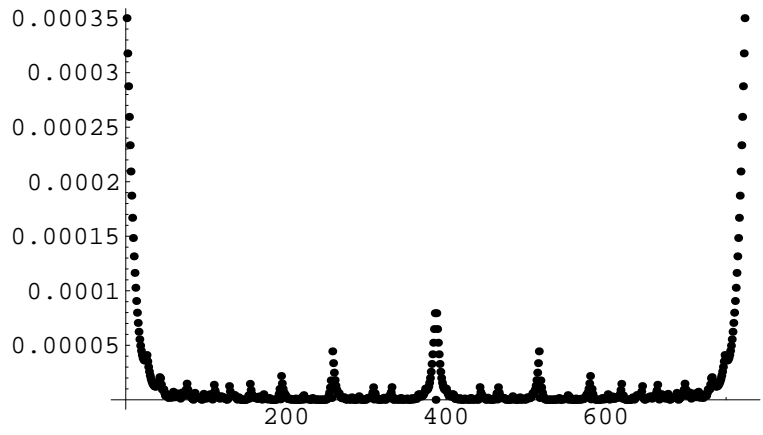}\\
  \end{center}\caption{Riemann-Farey Terms for $m = 50$}\label{FIG0}
\end{figure}

Figure~\ref{FIG0} is a plot of the terms in the sum (\ref{EQ1}) for
$m = 50$.

\pagebreak

Let $P_m(k)$ be sum of the $\phi(k)$ terms in (\ref{EQ1}) with Farey
denominator $k$ so that

\begin{equation}\label{EQ2}
R(m) = \sum_{k=2}^m P_m(k)
\end{equation}

\begin{figure}[h!]
\begin{center}
\leavevmode
  \includegraphics[width=100mm]{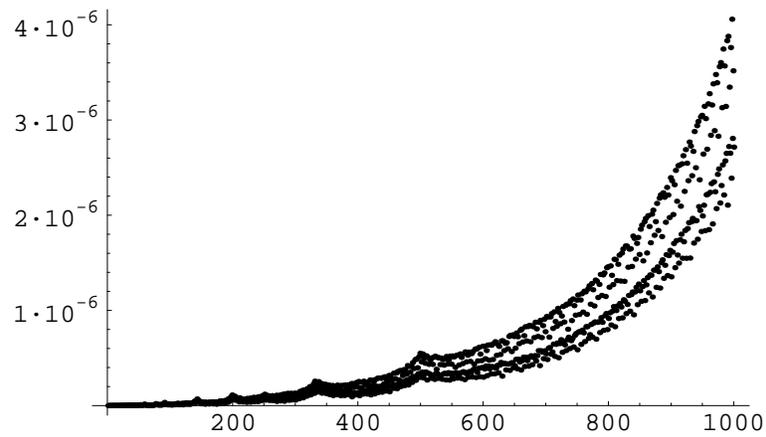}\\
  \end{center}\caption{$P_{m}(k)$ for $m=1000$}\label{FIG1}
\end{figure}

Figure~\ref{FIG1} is a plot of the terms in the sum (\ref{EQ2}) for
$m = 1000$.

The excursions from monotonicity (the ``bumps") in Figure~\ref{FIG1}
appear at abscissa values near $m/j$ for $2 < j < m$.
 \pagebreak

\begin{figure}[h!]
\begin{center}
\leavevmode
  \includegraphics[width=100mm]{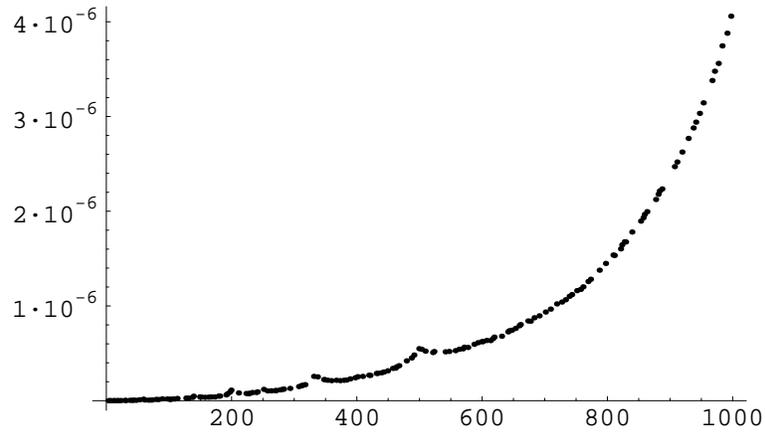}\\
  \end{center}\caption{$P_{m}(k)$ for $m=1000$ and prime $k$}\label{FIG1P}
\end{figure}

Figure~\ref{FIG1P} is the concave hull (the ``top") of
Figure~\ref{FIG1}. It is $P_m(k)$ for $k$ a prime.

We seek a closed form function of $k$, $\tilde{r}_m(k)$, that is a
tight upper bound on the points of Figure~\ref{FIG1P} so that

\begin{equation}
R(m) \leq \int_2^m {\tilde{r}_m(x)\: dx} \defas \tilde{R}(m)\text{.}
\end{equation}

\pagebreak

In general, for a fixed prime $m$ we compute $\tilde{r}_m(x)$ as

\begin{equation}
\tilde{r}_m(k) = e^{a_m + b_m k}\text{.}
\end{equation}

using values of $P_m(k)$.

Having done this for a collection of such $m$, we compute closed
forms for $a_m$ and $b_m$ as functions of $m$, $a(m)$ and $b(m)$
respectively, as

\begin{equation}
a(m) = s m^{t}
\end{equation}

and

\begin{equation}
b(m) = u m^{v}
\end{equation}

and study the asymptotic behavior of $R(m)$ by considering

\begin{equation}
\lim_{m\to\infty}R(m) \leq \lim_{m\to\infty}\tilde{R}(m) =
\lim_{m\to\infty} \int_2^m {e^{a(m) + b(m) m}}
\end{equation}

\pagebreak

In particular, for a fixed prime $m$ we compute $a_m$ and $b_m$
using the values of $P_m(k)$ at $k=m$ and at the prime $k$ closest
to $m/2$; i.e. at the right-most bump.

We then use the values of $a_m$ and $b_m$ obtained using values of
$m$ in sets of the form

\begin{equation}
M(p,q) = \left\{i^{th} \: \text{prime}\; |\; p \; \leq \; i \; \leq
\; q \right\}
\end{equation}

to compute values for $s$, $t$, $u$, and $v$ and hence $a(m)$ and
$b(m)$.

Table~\ref{TAB4} lists some $a(m)$ and $b(m)$ using this procedure.

\begin{table}[h!] \centerline{
\begin{tabular}{ccccc}
$m$ & $s$ & $t$ & $u$ & $v$\\
\hline
 &  &  &  &  \\
$M(101,200)$ & $-7.58$ & $0.112$ & $5.28$ & $-1.037$ \\
$M(201,300)$ & $-7.66$ & $0.111$ & $4.57$ & $-1.016$ \\
$M(301,400)$ & $-7.67$ & $0.111$ & $3.86$ & $-0.994$ \\
$M(401,500)$ & $-6.79$ & $0.125$ & $2.15$ & $-0.923$ \\
$M(501,600)$ & $-8.78$ & $0.094$ & $5.81$ & $-1.045$ \\
$M(601,700)$ & $-8.41$ & $0.099$ & $4.44$ & $-1.012$ \\
$M(701,800)$ & $-8.12$ & $0.103$ & $3.70$ & $-0.991$ \\
  &  &  &  &  \\
$M(101,800)$ & $-7.87$ & $0.107$ & $4.73$ & $-1.02$ \\
  &  &  &  &  \\
\hline
\end{tabular}}
\caption{Estimates of $a(m)$ and $b(m)$} \label{TAB4}
\end{table}

\pagebreak

\begin{figure}[h!]
\begin{center}
\leavevmode
  \includegraphics[width=100mm]{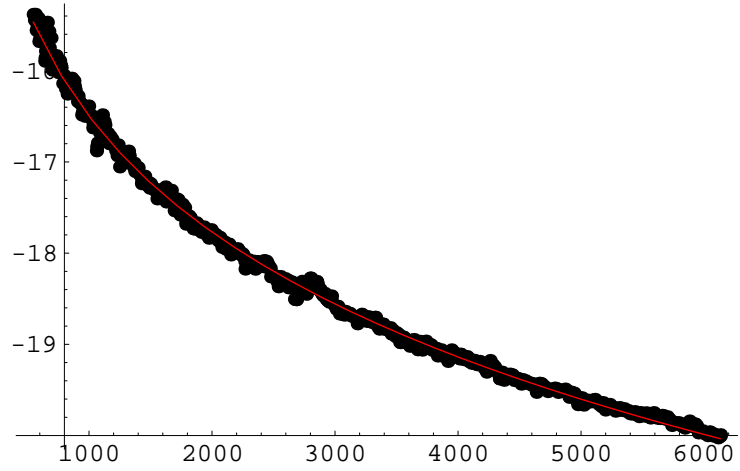}\\
  \end{center}\caption{$a(m)$ for $M(101,800)$ }\label{FIG8}
\end{figure}

\begin{figure}[h!]
\begin{center}
\leavevmode
  \includegraphics[width=100mm]{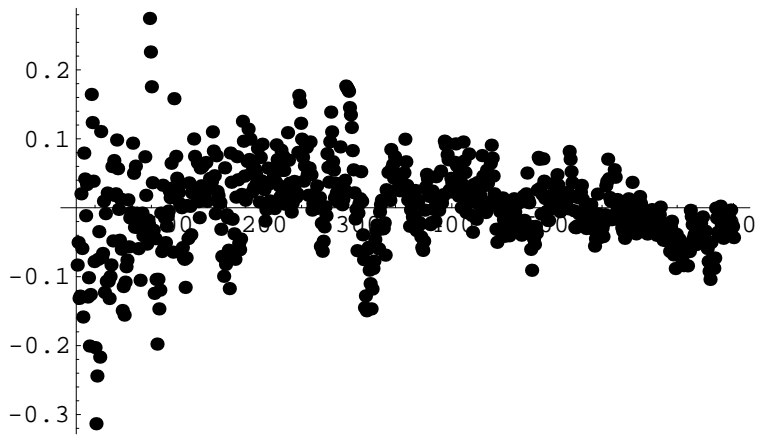}\\
  \end{center}\caption{Residuals of $a(m)$ for $M(101,800)$ }\label{FIG81}
\end{figure}

Figure~\ref{FIG8} is a plot of $a(m)$ for $M(101,800)$ and
Figure~\ref{FIG81} is a plot of the residuals.

\pagebreak

\begin{figure}[h!]
\begin{center}
\leavevmode
  \includegraphics[width=100mm]{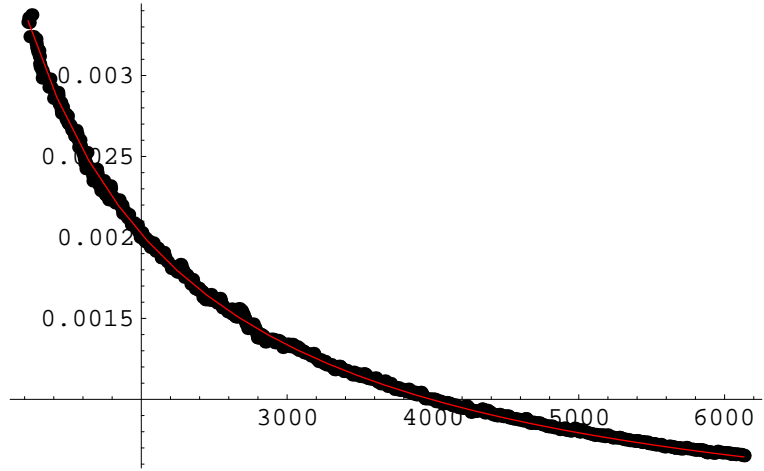}\\
  \end{center}\caption{$b(m)$ for $M(101,800)$}\label{FIG9}
\end{figure}

\begin{figure}[h!]
\begin{center}
\leavevmode
  \includegraphics[width=100mm]{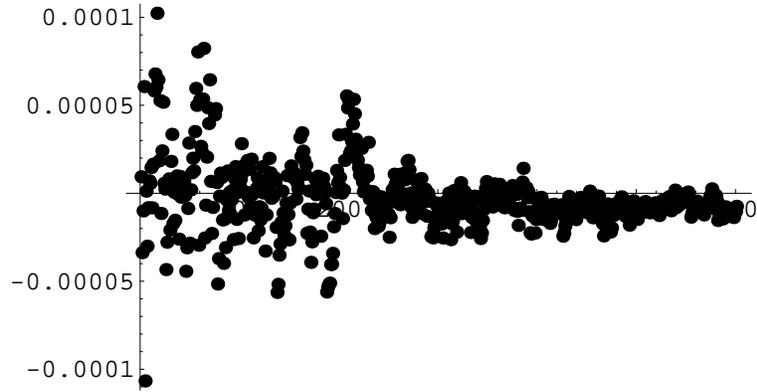}\\
  \end{center}\caption{Residuals of $b(m)$ for $M(101,800)$ }\label{FIG91}
\end{figure}

Figure~\ref{FIG9} is a plot of $b(m)$ for $M(101,800)$ and
Figure~\ref{FIG91} is a plot of the residuals.

\pagebreak

\begin{figure}[h!]
\begin{center}
\leavevmode
  \includegraphics[width=100mm]{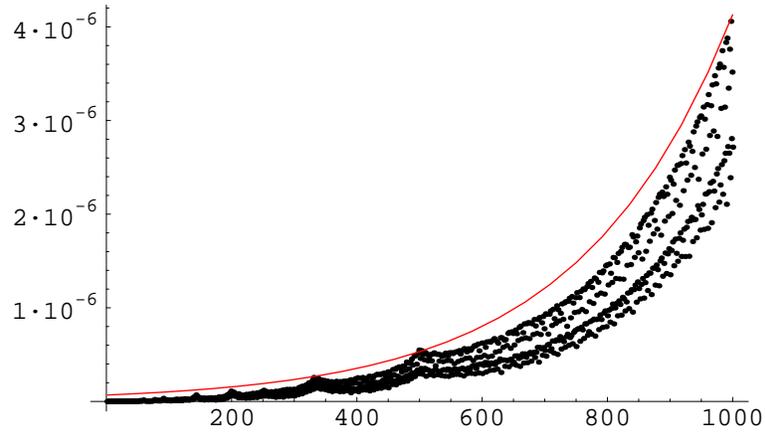}\\
  \end{center}\caption{$\tilde{r}_{m}(x)$ and $P_m(k)$ for $m=1000$}\label{FIG10}
\end{figure}

Taking $a(m)$ and $b(m)$ for $M(101,800)$, Figure~\ref{FIG10} is a
plot of

\begin{equation}
\tilde{r}_{m}(x) = e^{a(m) + b(m) x}
\end{equation}

and $P_m(k)$ for $m=1000$.

\pagebreak

Taking $a(m)$ and $b(m)$ for $M(101,800)$ we have the following
closed expression for $\tilde{R}(m)$:

\begin{equation}
\begin{split}
\tilde{R}(m) = \int_2^m {\tilde{r}_m(x) dx} = \int_2^m {e^{a(m) +
b(m) x} dx} = \\
-0.21 e^{-7.87 m^{0.11}}(e^{\alpha m^{1-\beta}}-e^{2 \alpha m
^{-\beta}})m^{\beta}
\end{split}
\end{equation}

where $\alpha = 4.73$ and $\beta = 1.02$.

\begin{figure}[h!]
\begin{center}
\leavevmode
  \includegraphics[width=100mm]{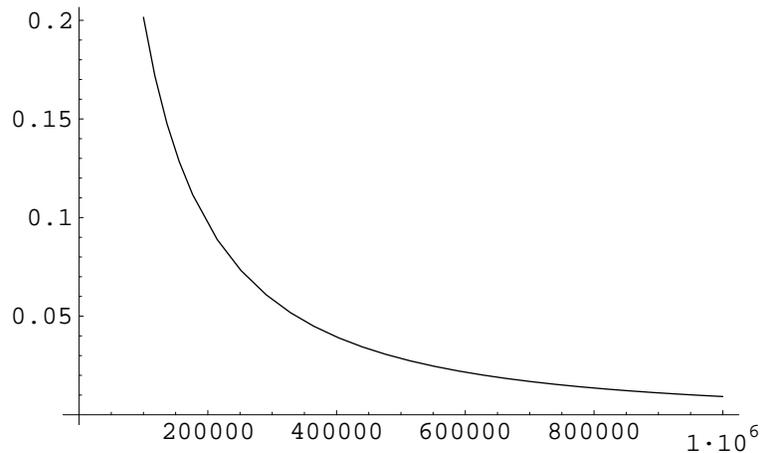}\\
  \end{center}\caption{$\hat R(x)/x^{-1+\epsilon}$}\label{FIG12}
\end{figure}

Figure~\ref{FIG12} is a plot of

\begin{equation}\label{EQ5}
\tilde R(x)/x^{-1+\epsilon}
\end{equation}

$x$ between $10^{5}$ and $10^{6}$ and $\epsilon = 0.000001$.

\pagebreak

To summarize, for the computations considered, we have

\begin{equation}
R(m)\leq \tilde R(m)
\end{equation}
and
\begin{equation}
\lim_{m\to \infty} \tilde R(m)/m^{-1+0.000001} = 0
\end{equation}

\end{document}